# New hysteresis operators with applications to counterterrorism

S. A. Belbas
Mathematics Department
University of Alabama
Tuscaloosa, AL. 35487-0350. USA.

e-mail: SBELBAS@GP.AS.UA.EDU, SBELBAS@BAMA.UA.EDU

Abstract.
We define two models of hysteresis that generalize the Preisach model. The first model is deterministic, the second model is stochastic and it utilizes discontinuous transition probabilities that satisfy impulsive differential equations. For the first model we prove, among other things, a local version of the "wiping out" property; for the stochastic model, we give methods for the construction of solutions of impulsive differential equations that determine the discontinuous transition probabilities. We also present a game-theoretic problem utilizing a generalized hysteresis operator. These hysteresis operators are motivated by questions of modelling the dynamics of decision making processes of networks of loosely knit terrorist groups.

1. Introduction.

This paper aims to provide a model of the decision-making processes of a loosely knit network of decision-making units. The underlying question is to model the behavior of a network of terrorist groups under conditions of "loose leadership". We can take a clue from other problems that involve situations where preferences of several units affect a collective result. Such situations have been studied, in the context of Economics, in references [CR, TMC] and other papers. The conclusion is that collective responses to stimuli can, in certain situations, be modelled by Preisach hysteresis operators. In general, a situation involving several agents, each of which makes a binary decision in response to an applied stimulus, according to some rules, naturally leads to hysteresis, in the following way: a binary decision process in response to a scalar-valued stimulus can be expressed as a non-ideal relay; for the purposes of the present paper, it is convenient to label the two possible decisions as $-1$ and $+1$ and describe a non-ideal relay in terms of two subsets $C_+$, $C_-$ of the real line $IR$, $C_- := (-\infty, \rho_1)$, $C_+ := (\rho_2, +\infty)$, $\rho_1 > \rho_2$, so that the decision changes from $-1$ to $+1$ when the stimulus exits from the set $C_-$, and switches from $+1$ to $-1$ when the stimulus exits from the set $C_+$; let us denote this relay by $R^{\rho_1 \rho_2}$ the total effect, as a resultant of the decisions of the individual units, can then be represented as an operator of the form $H = \iint\limits_{\rho_1 > \rho_2} w(\rho_1, \rho_2) R^{\rho_1 \rho_2} d\rho_1 d\rho_2$ where $w$ is a weighting function; this means that, if $u$ is an input signal, then

$$(Hu)(t) = \iint\limits_{\rho_1 > \rho_2} w(\rho_1, \rho_2)(R^{\rho_1 \rho_2} u)(t) d\rho_1 d\rho_2.$$

For the case of multi-dimensional input signals, and multiple, rather than binary, decisions, we need an extension of this model, and a corresponding analysis of the properties of such generalized models. A first step in this direction is undertaken in the next section.

The importance of modeling and analyzing the decision-making processes of networks of terrorist groups is well understood [C, CP, M, P, R, W]. Also, the works [G1, G2] utilize percolation studies to model the access of terrorists to targets via a social network; this last approach seems to us more suitable to homeland security issues rather than warfare in the terrorists' turf. Models of the decision processes and the behavior of terrorists can lead to techniques for data mining, in the sense of assimilating data into a theoretical model and identifying what kinds of data are useful and what kinds of missing data might lead to different analytical conclusions. In this paper, we aim at a mathematical analysis, and new mathematical constructs, based on our perception of what terrorists are, as gleaned from [AQ] and the factual information available in the news and in the quoted references. Our approach does not utilize existing models of terrorist decision making, but rather we look at the issue from a fresh point of view. Our main point of view is that terrorist groups have rational decision processes (rationality interpreted in a narrow sense), are connected through "loose leadership", and do not possess scientifically and computationally sophisticated systems of decision making. We model these characteristics by postulating that the decisions of each group, in response to external stimuli, can be represented by multi-dimensional multi-state relays, i.e. extensions of the one-dimensional two-state relay defined above, and that the cumulative effect of several groups, under conditions of "loose leadership", can be represented as an integral of a family of relays. These ideas are made mathematically precise in the next section.

## 2. Extension to vector signals and several output states.

When there are several output states, we denote the possible output states by lower-case Greek letters, and we denote by **A** the set of all elementary output states of a non-ideal relay. (In the case of the Preisach operator, there are two elementary output states; here, we consider an arbitrary number of elementary output states.) The input signal $u$ takes values in an open set $\Omega \subseteq IR^n$. A non-ideal relay is defined via a collection of open sets **C** := $\{C_\alpha : \alpha \in \mathbf{A}\}$, which we shall call <u>continuation sets</u>, with the properties

$$C_\alpha \subseteq \Omega \ \forall \alpha \in \mathbf{A}, \ \bigcup_{\alpha \in \mathbf{A}} C_\alpha = \Omega \tag{2.1}$$

and a partition of the relative boundary of each $C_\alpha$, i.e. the set $\mathrm{rel}\partial C_\alpha := \partial C_\alpha \cap \Omega$, into mutually disjoints subsets (some of which may be empty),

$$\mathrm{rel}\partial C_\alpha = \bigcup_{\beta \in \mathbf{A} \setminus \{\alpha\}} S_{\alpha\beta}, \ S_{\alpha\beta} \cap S_{\alpha\gamma} = \emptyset \text{ for } \beta \neq \gamma, \ S_{\alpha\beta} \subseteq C_\beta \ \forall \beta \in \mathbf{A} \setminus \{\alpha\}$$

$$\tag{2.2}$$

We set

$$\mathbf{S}_\alpha := \{S_{\alpha\beta} : \beta \in \mathbf{A} \setminus \{\alpha\}\} \text{ and } \mathbf{S} := \{\mathbf{S}_\alpha : \alpha \in \mathbf{A}\} \tag{2.3}$$

Given the above ingredients, the corresponding non-ideal relay $R$ is defined below. First, we need a definition: for every open set $C \subseteq \Omega$, for every signal $v$ taking values in $\Omega$, for every value t of the time-parameter, we define the <u>exit time</u> of $v$ from $C$ at times subsequent to $t$ by

$$\tau(t, v, C) := \inf\{s : s > t, v(s) \notin C\} \tag{2.4}$$

Clearly, this definition is nontrivial if and only if $v(t) \in C$.

For the signal $u$, taking values in $\Omega$ and defined for $t \in [t_0, T)$ where $T$ may be finite or $\infty$, we need to have an initial output state $\alpha_0 \in \mathbf{A}$ that satisfies the compatibility condition $u(t_0) \in C_{\alpha_0}$; we set $(Ru)(t_0) = \alpha_0$. Let $t_1 := \tau(t_0, u, C_{\alpha_0})$; if $t_1 \in I$, by our assumptions there exists a unique $\alpha_1 \in \mathbf{A} \setminus \{\alpha_0\}$ such that $u(t_1) \in C_{\alpha_1}$; then we define

$$(Ru)(t) = \alpha_0 \text{ for } t_0 \leq t < t_1;\ (Ru)(t_1^-) = \alpha_0;\ (Ru)(t_1^+) = \alpha_1 \tag{2.5}$$

This process is repeated inductively. For each $n$, we define $t_{n+1} := \tau(t_n, u, C_{\alpha_n})$; if $t_{n+1} \in I$, there exists a unique $\alpha_{n+1} \in \mathbf{A} \setminus \{\alpha_n\}$ such that $u(t_{n+1}) \in C_{\alpha_{n+1}}$; then we define

$$(Ru)(t) = \alpha_n \text{ for } t_n \leq t < t_{n+1}; \ (Ru)(t_{n+1}^-) = \alpha_n; \ (Ru)(t_{n+1}^+) = \alpha_{n+1} \qquad (2.6)$$

Once a non-ideal relay $R$ has been defined, we can formulate a superposition of non-ideal relays by analogy with the standard Preisach hysteresis operator. If $\rho$ is a parameter taking values in a measure space $(W, \mu)$ with $\mu(W)$ finite, we consider a family $\{(\mathbf{C}^\rho, \mathbf{S}^\rho): \rho \in W\}$ and the corresponding non-ideal relays $\{R^\rho: \rho \in W\}$. Given a signal $u$ and initial output states $\{\alpha_0^\rho: \rho \in W\}$ with the property that the functions $(R^\rho u)(t)$ are measurable as functions of $\rho$, we define the hysteresis operator $H$ by

$$(Hu)(t) := \int_W (R^\rho u)(t) \mu(d\rho) \qquad (2.7)$$

It should be noticed that, if the elementary output states (the elements of $\mathbf{A}$) and the measure $\mu$ are such that either the product $(R^\rho u)(t)\mu(d\rho)$ or the product $\mu(d\rho)(R^\rho u)(t)$ is vector valued (for example, for the second product, matrix-valued measure and vector-valued elementary output states), then the operator $H$ has vector-valued output.

The standard Preisach operator is clearly a particular case of the operator we defined in this section. The correspondence with the terminology and notation of this section is as follows: in the standard Preisach case, the dimension $n$ of the values of the input signals is $1$, and there are two elementary output states, which we may conventionally denote by $-1$ and $+1$; the sets $C_\alpha^\rho$ have the following form:

$C^\rho_{-1} = (-\infty, \rho_1)$, $C^\rho_{+1} = (\rho_2, +\infty)$, $\rho_1 \in IR$, $\rho_2 \in IR$, $\rho_1 > \rho_2$, $\rho = (\rho_1, \rho_2)$; the sets $S^\rho_{\alpha\beta}$ are

$S^\rho_{-1,+1} = \{\rho_1\}$, $S^\rho_{+1,-1} = \{\rho_2\}$; the set $W$ is $W = \{\rho = (\rho_1, \rho_2) \in IR^2 : \rho_L < \rho_2 < \rho_1 < \rho_U\}$

and $\mu$ is a finite Borel measure on $W$.

A natural question is: what are the essential properties of this hysteresis operator?

In the case of the standard Preisach operator, the following properties are known (see [M1]):

(i) the Preisach operator is causal and rate-independent;

(ii) the wiping-out property: the output of the Preisach operator is determined only by the history of dominant local maxima and dominant local minima of the input signal, and every other information is wiped out;

(iii) the property of congruence of first-order reversal loops;

(iv) there exists an algorithm, due to Mayergoyz [M1], for the identification of the weight function that appears in the definition of the Preisach operator;

(v) Mayergoyz's theorem: properties (i) through (iii) above are necessary and sufficient for an operator (acting on real-valued signals) to be representable as a Preisach hysteresis operator.

It is natural to want to know if there are counterparts of these properties for the operator defined in this section, and what modifications or new concepts are needed in order to formulate and prove such properties.

We have:

Theorem 2.1. Each relay $R$ is causal, i.e., if we denote by $(w)_t$ the restriction of any function $u$, defined on an interval $I$ that includes $t$, to the set $I \cap (-\infty, t]$, and if $u$ is an input signal with domain $D \equiv [t_0, T]$, then, for every $t$ in $D$ we have $(Ru)_t = (R(u)_t)_t$. Also, each relay $R$ is rate-independent, i.e., if $t = \varphi(s)$ is a change of the time variable, where $\varphi$ is a strictly increasing function, then $(Ru) \circ \varphi = R(u \circ \varphi)$.

Proof. The definition of $(Ru)(t)$ depends on the initial elementary output state of the relay and on the history, up to time $t$, of exit times of $u$ from the appropriate domains $C_\alpha$. For arbitrary but fixed $t$, these ingredients are the same for $u$ and for $(u)_t$. This proves the causality of $R$.

For the rate independence property, we observe that the exit times, say $\sigma(s, u \circ \varphi, C_\alpha)$, are related to the exit times $\tau(\varphi(s), u, C_\alpha)$ by $\tau(\varphi(s), u, C_\alpha) = \varphi(\sigma(s, u \circ \varphi, C_\alpha))$ since, by the strict monotonicity of $\varphi$, the condition $u(t') \notin C_\alpha$ is plainly equivalent to $(u \circ \varphi)(s') \notin C_\alpha$ & $t' = \varphi(s')$; consequently,

$$\tau(\varphi(s), u, C_\alpha) \equiv \inf\{t': t' \geq \varphi(s) \ \& \ u(t') \notin C_\alpha\} = \varphi(\inf\{s': s' > s \ \& \ u(\varphi(s')) \notin C_\alpha\}) \equiv$$
$$\equiv \varphi(\sigma(s, u \circ \varphi, C_\alpha).$$

///

Next, we formulate and prove an analogue to the wiping-out property. It turns out that our model has a property of <u>local wiping-out</u>. We first need some preliminary concepts and definitions. Let $U$ be an open subset of $\Omega$. Various local properties will refer to situations within $U$. On the set $W$ of indices $\rho$, we define the family of pre-orders $\leq_{(\alpha,U)}$ by

$\rho^{(1)} \leq_{(\alpha,U)} \rho^{(2)}$ if $C_\alpha^{\rho^{(1)}} \cap U \subseteq C_\alpha^{\rho^{(2)}} \cap U$. (Recall that a pre-order is reflexive and transitive, but not necessarily antisymmetric.)

An interval $(t_1, t_2)$ will be called an <u>interval of monotropy</u> for the signal $u$, relative to the pair of elementary output states $(\alpha, \beta)$, if, for all $t \in (t_1, t_2)$, the only possible switchings in $(R^\rho u)(t)$ are from $\alpha$ to $\beta$. A <u>maximal interval of monotropy</u> for $u$ relative to $(\alpha, \beta)$ is an interval of monotropy that is not properly contained in any other interval of monotropy. All time instants that do not belong to an interval of monotropy are called <u>transition points</u>.

We have:

<u>Proposition 2.1.</u> For each $x$ in $\Omega$, let $d(x) := \inf_{\rho, \alpha, \beta} \text{dist}(x, S_{\alpha\beta}^\rho)$, where the infimum is taken over all $\rho$ in $W$ and all $\alpha$, $\beta$ in **A** with $\alpha \neq \beta$. Let $u$ be a continuous input signal taking values in $\Omega$. If, at some time-instant $t'$, where $t'$ is an interior point of the domain of $u$, we have $d(u(t'))>0$, then $t'$ belongs to an interval of monotropy.

Proof. Let $d':=d(u(t'))>0$. For every $d''<d'$, we have $B(u(t'),d'') \cap S^\rho_{\alpha\beta} = \emptyset$ for all $\alpha, \beta, \rho$, where $B(x,r) := \{y \in \mathbb{R}^n : |y-x| < r\}$. Since $\Omega$ is open, there is a $d'''$, $0 < d''' < d'$, such that $B(u(t'),d''') \subseteq \Omega$. By the continuity of $u$, there is an $\varepsilon > 0$ such that $u(t) \in B(u(t'),d''') \ \forall t \in (t'-\varepsilon, t'+\varepsilon)$, thus, for every $\rho, \alpha, \beta$, we have $u(t) \notin S^\rho_{\alpha\beta} \ \forall t \in (t'-\varepsilon, t'+\varepsilon)$. According to our definitions above, $t'$ belongs to an interval of monotropy for $u$, for every $\alpha$ and $\beta$. ///

Now, we assume that $U$ is an open set in $\Omega$ with the following properties:

(i) $U \cap \mathbf{C}$ (by which we mean the family $\{U \cap C^\rho_\alpha : \alpha \in \mathbf{A}, \rho \in W\}$) consists only of sets of the form $U \cap C^\rho_{\alpha_0}$, $U \cap C^\rho_{\alpha_1}$, $\rho \in W' \subseteq W$ for two arbitrary but fixed indices $\alpha_0, \alpha_1$ in $\mathbf{A}$.

(ii) The hysteresis operator $H$, acting on an input signal that takes values in $U$, evolves according to the following rules:

at an initial time $t_0$, all relays are at the elementary output state $\alpha_0$; as soon as $u(t)$ reaches a relative boundary $S^\rho_{\alpha_0} \cap U$, for some $\rho$ in $W'$, the relay $R^\rho$ is turned to the elementary output state $\alpha_1$; continuing in this way, every time $u(t)$ reaches a relative boundary $S^\rho_{\alpha_i} \cap U$, say at time $t'$, all relays $R^\rho$ that satisfy $(R^\rho u)(t'^-) = \alpha_i$, where $(R^\rho u)(t'^-) \equiv \lim_{t \to t'^-} (R^\rho u)(t)$, will be switched to state $\alpha_{i+1 (\mathrm{mod}\, 2)}$.

(iii) Whenever a relay $R^\rho$ is witched from state $\alpha_i$ to

$\alpha_j$, all relays $R^{\rho'}$, with $\rho' \in W'$ & $\rho' <_{(\alpha_i, U)} \rho$, will be at state $\alpha_j$.

We have the following:

Theorem 2.2. (Local wiping-out property.) Let $U$ be an open subset of $\Omega$, and $u$ be an input signal, defined for $t \in [t_0, T)$, such that $u(t)$ lies in $U$ for all $t \in [t_0, T)$. We assume that conditions (i) through (iii) above are satisfied, and that $u$ has a finite number of isolated transition points that alternate between transition points for $(\alpha_0, \alpha_1)$ and transition points for $(\alpha_1, \alpha_0)$. Then the values of $(Hu)(t')$, at each transition point $t'$, depend only on the history of transition points up to and including time $t'$. If $t'$, $t''$ are two transition points, relative to $S^{\rho'}_{\alpha_i \alpha_{i+1 (\mathrm{mod}\, 2)}}$, $S^{\rho''}_{\alpha_i \alpha_{i+1 (\mathrm{mod}\, 2)}}$, respectively, and if

$t'' > t'$ & $\rho'' >_{(\alpha_i, U)} \rho'$, then the value of $(Hu)(t'')$ depends on the history of transition points up to and including $t''$ but excluding $t'$, i.e. the effect of the transition point $t'$ has been wiped out at time $t''$.

Proof. We examine what happens between two consecutive transition points. Let $t_1$, $t_2$ be two consecutive transition points, relative to $S^{\rho^{(1)}}_{\alpha_0\alpha_1}$, $S^{\rho^{(2)}}_{\alpha_1\alpha_0}$. At time $t_1$, the state of $(R^{\rho^{(1)}}u)$ is changed from $\alpha_0$ to $\alpha_1$. For $t_1 < t < t_2$, the only possible switchings in $(R^\rho u)$ are from $\alpha_1$ to $\alpha_0$, since, by assumption, there is no transition point between $t_1$ and $t_2$. For $t_1 < t < t_2$, whenever $u(t) \in S^\rho_{\alpha_1\alpha_0}$, all relays $R^{\tilde\rho}$ with $\tilde\rho \leq_{(\alpha_1, U)} \rho$ that were previously at state $\alpha_1$ have been turned to state $\alpha_0$. By the time $t_2$, all relays $R^{\tilde\rho}$ with $\tilde\rho <_{(\alpha_1, U)} \rho^{(2)}$ that were previously at state $\alpha_1$ have been turned to state $\alpha_0$, whereas the relay $R^{\rho^{(2)}}$ will be switched from $\alpha_1$ to $\alpha_0$ precisely at time $t_2$. Thus the relays that have been switched between times $t_1$ and $t_2$ depend only on two things: first, which relays were in each of the two states at time $t_1$, and second the position of $u(t_2)$, of course taking into account the information that $(t_1, t_2)$ is a maximal interval of monotropy.

Now, consider two transition points, say $t'$ and $t''$, both transition points for the same ordered pair, say $(\alpha_0, \alpha_1)$. According to our assumptions, $t'$ and $t''$ are not consecutive transition points. For simplicity, we treat the case when there is just one more transition point $t^\wedge$ in $(t', t'')$, and $t^\wedge$ is a transition point relative to $(\alpha_1, \alpha_0)$; other cases can be treated similarly. Let $\rho'$, $\rho''$, $\rho^\wedge$ be the values of $\rho$ corresponding to the transition points $t'$, $t''$, $t^\wedge$. By assumption, $\rho' <_{(\alpha_0, U)} \rho''$. At time $t'$, all relays $R^\rho$ with $\rho \leq_{(\alpha_0, U)} \rho'$ that were previously at state $\alpha_0$ are switched to the state $\alpha_1$; an unspecified set of relays will be in state $\alpha_0$. The relays that had been turned to either state at the transition point before $t'$ are treated as a fixed set of relays, say **R**. Thus the relays that will be switched from state $\alpha_0$ to $\alpha_1$ at time $t'$, out of the set **R**, are precisely those that satisfy $\rho \leq_{(\alpha_0, U)} \rho'$. At time $\hat{t}$, those of the relays that were at state $\alpha_1$ and satisfy $\rho \leq_{(\alpha_1, U)} \hat{\rho}$ will be switched to state $\alpha_0$. At time $t''$, out of the relays that were at state $\alpha_0$, those that satisfy $\rho \leq_{(\alpha_0, U)} \rho''$ (and that set includes those relays that had been switched from $\alpha_0$ to $\alpha_1$ at time $t'$, since $\rho' <_{(\alpha_0, U)} \rho''$ and the switching at time $\hat{t}$ introduced possibly more relays in the $\alpha_0$ state) will be switched from $\alpha_0$ to $\alpha_1$. Those relays, out of the relays that were switched from $\alpha_1$ to $\alpha_0$ at time $\hat{t}$ with the $\alpha_1$ not due to the switchings at $t'$, are those that would have been switched from $\alpha_1$ to $\alpha_0$ without the switchings at time $t'$; those relays, out of the relays that were switched from $\alpha_1$ to $\alpha_0$ at time with the $\alpha_1$ coming from the switchings at $t'$, will be switched back to state $\alpha_1$ at time $t''$. Thus the output $(Hu)(t'')$ depends on the situation at the transition point before $t'$, the transition point $\hat{t}$, and the situation at $t''$. In this way, the effect of the switchings at time $t'$ has been wiped out at time $t''$. ///

In the above theorem and proof, the transition points play a role analogous to local extrema of the input signal in the case of the standard Preisach model. The condition $\rho' <_{(\alpha_i, U)} \rho''$ is related to the concept of a dominant local extremum of the standard Preisach operator.

We give an example of how a decision rule of the form of a relay with more than two states might be formulated. Obviously, this is a hypothetical example, intended to show just one possible concrete instance of the model we have presented.

Example 2.1. We consider a two-dimensional input signal $u = (u^1, u^2)$ satisfying the conditions

$$u^1 > 0, \; u^2 > 0, \; a^1 u^1 + a^2 u^2 < c \quad (a^1 > 0, \; a^2 > 0, \; c > 0)$$

The inputs $u^1$ and $u^2$ represent two types of actions, economic sanctions ($u^1$) and military operations ($u^2$). The domain $\Omega$ on the $u^1 u^2$-plane is the interior of the triangle $A_1 A_2 A_3$ in figure 2.1.

**Figure 2.1 goes here.**

The relay has 3 elementary output states, representing 3 levels of actions on the part of the terrorist group; for instance, if there are two possible types targets, the 3 possible actions could be $\alpha_1 :=$ (hit both types of targets, #1 and #2), $\alpha_2 :=$ (hit targets of type #2 only), $\alpha_3 :=$ (hit targets of type #1 only). The threshold boundaries $S_i$, $i = 1,2,3$, are given by

$$S_i = \{(u^1, u^2) \in \Omega : a_i^1 u^1 + a_i^2 u^2 = c_i\}$$

We shall use the notation $(PQ)$ to denote the non-oriented (i.e. $(PQ) \equiv (QP)$) open straight line segment with endpoints $P$ and $Q$, $[PQ]$ for the corresponding non-oriented closed line segment, and $[PQ)$, $(PQ]$ for the non-oriented semi-open segments. The boundaries $S_1$, $S_2$, $S_3$ are the line segments $(A_{23}A_{32})$, $(A_{31}A_{13})$, $(A_{12}A_{21})$, respectively, in figure 2.1. The sets $C_1$, $C_2$, $C_3$ are the interiors of the convex quadrilaterals $A_2A_3A_{32}A_{23}$, $A_3A_1A_{13}A_{31}$, $A_1A_2A_{21}A_{12}$, respectively, in the same figure. Naturally, in order to have a meaningful relay, the relative positions and slopes of the segments $(A_{23}A_{32})$, $(A_{31}A_{13})$, $(A_{12}A_{21})$ have to be such that the covering condition $C_1 \cup C_2 \cup C_3 = \Omega$ is satisfied, as in the figure. The points $B_1$, $B_2$, $B_3$ are intersections of the pairs of segments

$((A_{12}A_{21}), (A_{31}A_{13}))$, $((A_{23}A_{32}), (A_{12}A_{21}))$, $((A_{31}A_{13}), (A_{23}A_{32}))$, respectively.

We assume that the parameters $a_i^j$, $c_i$, $i = 1,2,3$, $j = 1,2$, are positive, and therefore the continuation sets $C_i$ are given by $C_i = \{u \in \Omega : a_i^1 u^1 + a_i^2 u^2 < c_i\}$ for $i=2, 3$, and

$$C_1 = \{u \in \Omega : a_1^1 u^1 + a_1^2 u^2 > c^1\}.$$

The situation described above has the following implications: as long as the combinations $a_i^1 u^1 + a_i^2 u^2$, $i=2, 3$ do not exceed the values $c_i$, the terrorists will attack only one type of targets; the choice of type of targets depends on the values of $u$ and the history of previous decisions of the terrorists; when the value of $a_1^1 u^1 + a_1^2 u^2$ exceeds $c_1$, the terrorists will attack both types of targets; for other situations, choices are made according to the rules set forth in our general definition of a relay. It will be noticed that the model is not yet complete, because there are parts of the boundaries where the choice is still ambiguous, namely the segements $(B_1B_2)$, $(B_2B_3)$, $(B_3B_1)$, each of which belongs to the intersection of two continuation sets. In order to complete the model, we take points $D_1$, $D_2$, $D_3$ on the segnents $(B_2B_3)$, $(B_3B_1)$, $(B_1B_2)$, respectively, and we take

$$S_{12} := [D_1A_{32}),\ S_{13} := (D_1A_{23});$$
$$S_{31} := [D_2A_{21}),\ S_{32} := (D_2A_{12});$$
$$S_{23} := [D_3A_{13}),\ S_{21} := (D_3A_{31})$$

The points $D_1$, $D_2$, $D_3$ are chosen in order to remove ambiguity (in the model) about switching actions when the input signal, at the moment of exit from one continuation set, happens to lie in the intersection of the two other continuation sets. In the actual behavior of the terrorist group, that ambiguity may be non-removable; in that case, the area of the triangle $B_1B_2B_3$ gives a measurable assessment of the level of ambiguity. Now we have defined a relay operator that describes the decision-making processes of the terrorist as their response to the input signal *u*.  ///

## 3. Extension to Markov processes.

The model of the previous section can be extended to a stochastic framework, by using the idea of changing sets into Markov transition probabilities.

To motivate this idea, consider first a modified version of the situation described in the previous section. In this modification, whenever a part of the boundary $S_\alpha$ belongs to more than one set $C_\beta$, the choice of switching action is not deterministic, but rather, at each point $x$ of $S_\alpha$, there is a probability $p_{\alpha\beta}(x)$, with $\beta$ running over all values for which $x \in C_\beta$, of switching to elementary output state $\beta$. Immediately after the switching, the elementary output state of the relay is not a deterministic state, but becomes a probability distribution; at the next step, we are having a stochastic output state, in other words, a random variable taking values into a set of various elementary output states, with a given probability distribution; at the next step, for each of those output states, by the time the input reaches a new boundary $S_\beta$, there will be a new set of transition probabilities, for the stochastic switching from $\beta$ to a new elementary output state. Thus, in this scenario, the role of the boundaries in **S** is to define discontinuities (or jumps) in the probability distribution of the stochastic elementary output state.

For this reason, a consistent and relatively simple formulation can be achieved by combining all the different boundaries in **S** from the model into one single set $S$, and then specify the nature of the discontinuities in the probability distribution of the elementary output state of a relay.

This is analogous to the idea used in randomized testing of statistical hypotheses, where the rejection set (points of sample observations that lead to rejection of the null hypothesis) is replaced by a probability distribution. A set $C_\alpha$ can be identified with its own indicator function $\chi_{C_\alpha}$, which is defined as $\chi_{C_\alpha}(x) := \begin{cases} 1, & \text{if } x \in C_\alpha \\ 0, & \text{if } x \notin C_\alpha \end{cases}$.

Similarly, every set $S_{\alpha\beta}$ can be identified with the indicator function $\chi_{S_{\alpha\beta}}$. Because of the assumption $S_{\alpha\beta} \subseteq C_\beta$, the indicator function $\chi_{S_{\alpha\beta}}$ is the restriction to $S_{\alpha\beta}$ of the function $\chi_{C_\beta}$.

The analogue of randomization, in our case, is to replace $\chi_{C_\alpha}$ and $\chi_{S_{\alpha\beta}}$ by functions $p_{\alpha\alpha}(t,x)$, $p_{\alpha\beta}(t,x)$, $\beta \neq \alpha$, which satisfy the nonnegativity and stochasticity conditions

$$p_{\alpha\beta}(t,x) \geq 0, \sum_\beta p_{\alpha\beta}(t,x) = 1.$$

These functions are defined for $(t, x) \in S \subseteq (0,T) \times \Omega$.

Given an input $u(t)$, the output of this relay is a matrix of Markov transition probabilities. These transition probabilities depend on spatial location $x$ in $\Omega$, and they will be denoted by $\pi_{\alpha\beta}(s,t,x)$. For $(t,x) \notin S$, the family $\pi_{\alpha\beta}(s,t,x)$ satisfies the standard (forward) Kolmogorov equations (cf., e.g., [BR])

$$\frac{\partial \pi_{\alpha\beta}(s,t,x)}{\partial t} = -\varphi_{\beta}(t,x)\pi_{\alpha\beta}(s,t,x) + \sum_{\gamma \in \mathbf{A}} \varphi_{\gamma}(t,x) g_{\gamma\beta}(t,x) \pi_{\alpha\gamma}(s,t,x) \; ;$$

$$\pi_{\alpha\beta}(s,s,x) = \delta_{\alpha\beta} \equiv \begin{cases} 1, & \text{if } \alpha = \beta \\ 0, & \text{if } \alpha \neq \beta \end{cases}$$

(3.1)

The various terms above are interpreted in the context of the theory of Markov processes in continuous time with discrete and finite state space. The set **A** is the state space of a Markov process. The functions $\varphi_{\alpha}(t,x)$ are the so-called intensities: the probability hat the state (value) of the process will change in the interval $(t, t+\delta t)$, conditioned on the event that the process is at state $\alpha$ at time t, is $\varphi_{\alpha}(t,x)\delta t + o(\delta t)$. The functions $g_{\alpha\beta}(t,x)$ represent the conditional probabilities that the process will be at state $\beta$ at time $t+\delta t$, conditioned on the events that the process is at state $\alpha$ at time $t$ and that the state of the process changes in the interval $(t, t+\delta t)$. At points $(t,x) \in S$, the transition probabilities $\pi_{\alpha\beta}$ change according to the impulse condition

$$\pi_{\alpha\beta}(s,t^{+},x) = \sum_{\gamma \in \mathbf{A}} \pi_{\alpha\gamma}(s,t^{-},x) p_{\gamma\beta}(t,x) \tag{3.2}$$

Thus, the time-evolution of the transition probabilities is given by a system of impulsive ordinary differential equations. An exposition of the main results and techniques of impulsive differential equations may be found in [BS].

A further modification in the Kolmogorov equations is what we shall call <u>the semi-flow version</u>. The variable x is replaced by a continuous semi-flow $u(s,t;\xi)$, $0 \leq s \leq t \leq T$, $\xi \in IR^n$. The defining properties of a semi-flow are

$$u(s,s;\xi) = \xi \text{ and } u(s,t;u(r,s;\xi)) = u(r,t;\xi) \text{ for } r \leq s \leq t \quad (3.3)$$

The semi-flow version of the Kolmogorov equations with impulses is

$$\frac{\partial \pi_{\alpha\beta}(s,t,u(s,t;\xi))}{\partial t} = -\varphi_\beta(t,u(s,t;\xi))\pi_{\alpha\beta}(s,t,u(s,t;\xi)) +$$

$$+ \sum_{\gamma \in \mathbf{A}} \varphi_\gamma(t,u(s,t;\xi))g_{\gamma\beta}(t,u(s,t;\xi))\pi_{\alpha\gamma}(s,t,u(s,t;\xi)) ; \quad (3.4)$$

$$\pi_{\alpha\beta}(s,s,\xi) = \delta_{\alpha\beta} ;$$

$$\pi_{\alpha\beta}(s,t^+,u(s,t;\xi)) = \sum_{\gamma \in \mathbf{A}} \pi_{\alpha\gamma}(s,t^-,u(s,t;\xi))p_{\gamma\beta}(t,u(s,t;\xi)) \quad (3.5)$$

Clearly, the model (3.4, 3.5) ) is a particular case of (3.1, 3.2)) above.

The role of a semi-flow is to represent the input signal to this stochastic version of a non-ideal relay.

Given an input $u(t)$, $0<t<T$, with values in $\Omega$, the output of this stochastic relay is

$$(Ru)(s,t) = \{\pi_{\alpha\beta}(s,t,u(s,t;\xi)): \alpha, \beta \text{ in } \mathbf{A}\} \tag{3.6}$$

Given a family of stochastic relays $R^\rho$, $\rho \in W$, with corresponding transition probabilities $\pi^\rho(s,t,u(s,t;\xi))$, the <u>stochastic hysteresis operator</u> is defined as

$$(Hu)(s,t) = \{\int_{\rho \in W} \pi^\rho(s,t,u(s,t;\xi))\mu(d\rho): \alpha, \beta \text{ in } \mathbf{A}\} \tag{3.7}$$

To the best of our knowledge, the concepts of transition probabilities that satisfy impulsive variants of the Kolmogorov equations, and semi-flow extension of these equations, are new concepts introduced here for the first time.

Next, we describe methods for solving the impulsive variant of the Kolmogorov equations. First, we simplify the notation; we suppress the dependence of the various terms on $s$, $u$, and $\xi$, and we write $\psi(t)$ for the transpose of the matrix $\pi(t)$. We assume that, for each $s$, $u$, and $\xi$, there is a finite set $\{\tau_1, \tau_2, ..., \tau_N\}$ of time instants $\tau$ in $(s, T)$ that satisfy $(\tau, u(s,\tau;\xi)) \in S$. This means that the curve $\{(\tau, u(s,\tau;\xi)): s < \tau < T\}$ intersects the set *S* at a finite number of points. Then the impulsive Kolmogorov equations have the form

$$\frac{d\psi(t)}{dt} = A(t)\psi(t) \text{ for } t \notin \{\tau_1,...,\tau_N\};$$
$$\psi(t^+) = B(t)\psi(t^-) \text{ for } t \in \{\tau_1,...,\tau_N\}; \quad (3.8)$$
$$\psi(s) = I \ (I = \text{identity matrix})$$

For notational convenience, we set $\tau_0 := s$, $\tau_{N+1} := T$. The construction of a fundamental matrix for impulsive systems of the type of (3.8) has been carried out in [BS], and we outline here this method, with suitable notational changes.

For each *t'*, *t* in $(\tau_i, \tau_{i+1}]$, *i=0, 1, 2, ..., N*, with *t'<t*, we denote by $\Psi(t', t)$ the fundamental matrix of the system $\frac{d\psi(t)}{dt} = A(t)\psi(t)$, i.e. $\Psi(t', t)=I$ and

$$\frac{\partial \Psi(t',t)}{\partial t} = A(t)\Psi(t',t) \ \forall t \in (\tau_i, \tau_{i+1}].$$ Also, we define $\Psi(\tau_i^+, t) := \lim_{t' \to \tau_i^+} \Psi(t', t)$. Then,

for $\tau_{i-1} < t' \leq \tau_i < \tau_j < t \leq \tau_{j+1}$ the fundamental matrix of the impulsive system is

$$\Phi((t')^+, t) = \Psi((t')^+, \tau_i) \left[ \prod_{l=i}^{j-1} \left( B(\tau_l) \Psi(\tau_l^+, \tau_{l+1}) \right) \right] \Psi(\tau_j^+, t)$$

This construction depends on knowing the fundamental matrices $\Psi(t', t)$ on every interval $(\tau_i, \tau_{i+1}]$. A constructive method, that does not require knowledge of $\Psi(t', t)$, but produces a sequence of successive approximations to the fundamental matrix $\Phi(t', t)$, is available as a particular case of the technique introduced in [BS1] for a related problem for Volterra integral equations with impulses. The specialization of the method of [BS1] to problem (3.8) above amounts to the following: for $j>i$, we define the set $\mathbf{P}(i,j)$ of increasing paths from $i$ to $j$ as the set of all collections of indices $i < k_1 < k_2 < ... < k_r < j$; for $j=i$, we define $\mathbf{P}(i,i):=\{i\}$; for every $\sigma$ in $\mathbf{P}(i,j)$, we define $V(\sigma)$ by

$$V(\sigma) := B(\tau_{k_r})B(\tau_{k_{r-1}})...B(\tau_{k_1})B(\tau_i) \text{ if } i < j \text{ and } \sigma = \{i, k_1, ...., k_r, j\};$$
$$V(\sigma) := 1 \text{ if } \sigma = \{i\};$$

we define the kernels $K(t,t')$ by

$$K(t,t') := A(t') + \sum_{j:s<\tau_j<t} \sum_{i=1}^{j} \sum_{\sigma \in P(i,j)} B(\tau_j)V(\sigma)A(t);$$

we define the convolution of any two kernels having the form of K(t, t') above, by

$$(K_1 * K_2)(t,t') := \int_{t'}^{t} K_1(t,t_1)K_2(t_1,t')dt_1;$$

then $\Phi(t',t) = I + \sum_{m=1}^{\infty} \int_{t'}^{t} K^{*m}(t,t_1)dt_1$, where the exponent "$*m$" denotes $m$ convolutions of $K$ with itself.

## 4. Combining hysteresis with game theory.

The purpose of this section is to illustrate the possibility of utilizing the deterministic generalized hysteresis operator, introduced in section 2, in the formulation and analysis of differential games. It must be emphasized that the formulation below is intended as a demonstration of the possibility of using generalized hysteresis operators in a game-theoretic situation, and it is not represented as an example that can be directly applied to an actual real-world situation involving terrorist warfare; such an example would require much more work, with information and data that are not presently available to the author.

We consider a differential game with two players, $P_1$ and $P_2$, with corresponding controls $u_1(t), u_2(t), 0 \leq t \leq T$. Player $P_2$ corresponds to a collection of terrorist groups under conditions of "loose leadership". The control $u_2$ consists of two parts, $u_2(t) = (u_{21}(t), u_{22}(t))$, where $u_{21}(t) = g(t, (Hu_1)(t))$, $H$ being a hysteresis operator of the type defined in section 2, and $g$ a (generally nonlinear) function. This representation of $u_2$ is intended to describe a situation where some of the decisions are made on the basis of rational calculations and coordination of the constituent groups (these decisions are represented by $u_{22}$), whereas others are made under conditions of "loose leadership" and/or rigid criteria, and they can be represented as nonlinear functions of a hysteresis operator, acting on the signal $u_1$ which is the control decided by player $P_1$. The state $y(t)$ of the system that is being controlled by both players satisfies

$$\frac{dy(t)}{dt} = f(t, y(t), u_1(t), u_2(t)); \ y(0) = y_0 \tag{4.1}$$

The performance functional $J$ is given by

$$J := \int_0^T F(t, y(t), u_1(t), u_2(t))dt + F_0(y(T)) \qquad (4.2)$$

and we seek $\sup_{u_1} \inf_{u_{22}} J$.

For simplicity, we assume that $W$ is a finite set, and we use the notation $\overline{\alpha} := (\alpha_\rho : \rho \in W)$ for elements in $\mathbf{A}^W$. We denote by $V^{\overline{\alpha}}(t,x)$ the parameterized value function of the game, i.e. the value of the corresponding game starting at time $t$ in state $x$ while the relays are in state $\overline{\alpha}$, i.e. each relay $R^\rho$ is in state $\alpha_\rho$. We set

$$S^{\overline{\alpha}\overline{\beta}} := \bigcap_{\rho \in W} S_{\alpha_\rho \beta_\rho}, \quad S^{\overline{\alpha}} := \bigcup_{\overline{\beta}:\overline{\beta} \neq \overline{\alpha}} S^{\overline{\alpha}\overline{\beta}} \qquad (4.3)$$

The expression $g(t, (Hu_1)(t))$ is clearly a function of $t$ and $\overline{\alpha}$, and we set $G(t, \overline{\alpha}) \equiv g(t, (Hu_1)(t))$. We denote by $c_1$, $c_2$ the admissible values of $u_1$, $u_{22}$, respectively. Then, formally, the value function $V$ satisfies

$$\sup_{c_1} \inf_{c_2} \{ \frac{\partial \overline{V^\alpha}(t,x)}{\partial t} + \sum_i \frac{\partial \overline{V^\alpha}(t,x)}{\partial x_i} f_i(t,x,c_1,(G(t,\overline{\alpha}),c_2)) +$$

$$+ F(t,x,c_1,(G(t,\overline{\alpha}),c_2)) \} = 0 \text{ for } c_1 \notin S^{\overline{\alpha}},\ 0 \leq t < T; \quad (4.4)$$

$$\overline{V^\alpha}(t^-,x) = \overline{V^\beta}(t^+,x) \text{ for } c_1 \in S^{\overline{\alpha}\overline{\beta}};\ \overline{V^\alpha}(T,x) = F_0(x)$$

This equation is, in some respects, the game-theoretic counterpart of the dynamic programming equations for optimal control problems for systems with hysteresis [BM1, BM2, BM3]. The analysis and solution methodology for game-theory equations (which are analogous to the Bellman-Isaacs equation of ordinary differential games) is a topic for further research.

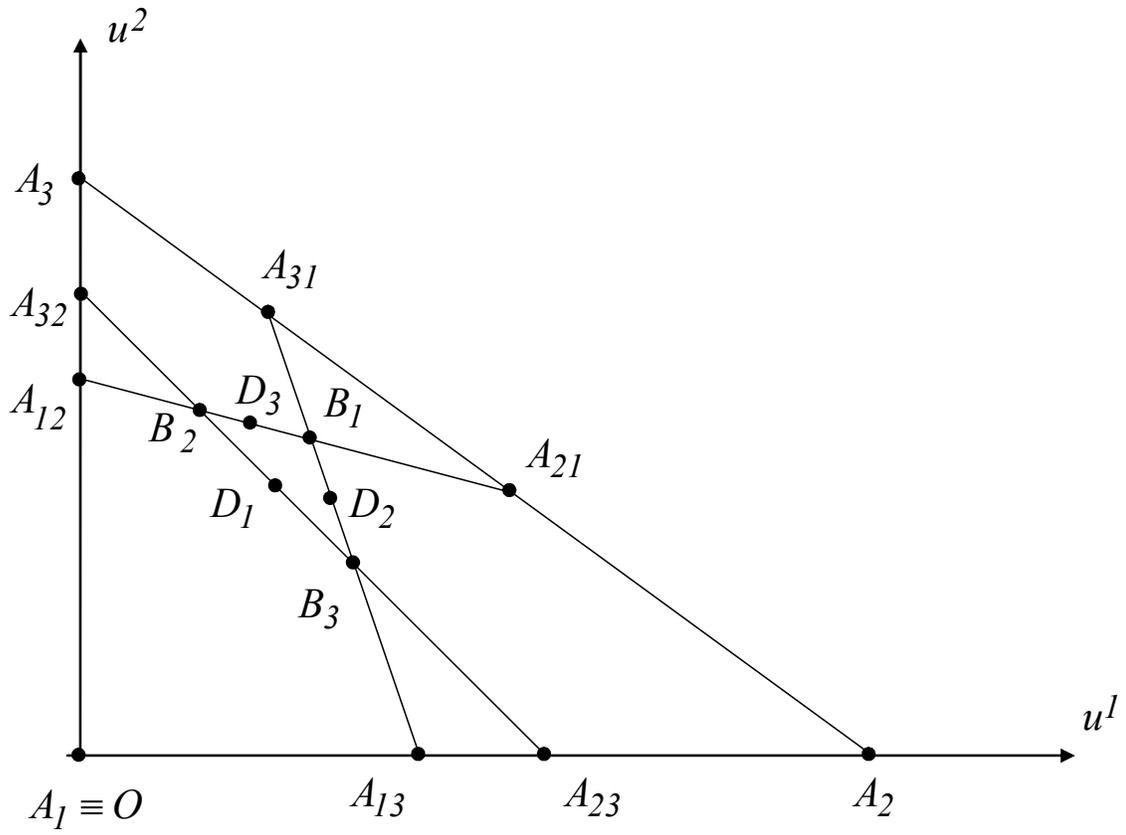

FIGURE 2.1.